\documentclass[10pt,draft]{article}
\usepackage{amsfonts,amsmath,amssymb,cite}
\newcommand{\Li}{\mathop{\textrm{Li}}\nolimits}

\begin{document}
\title{The parameter derivatives 
$[\partial^{2}P_{\nu}(z)/\partial\nu^{2}]_{\nu=0}$ and 
$[\partial^{3}P_{\nu}(z)/\partial\nu^{3}]_{\nu=0}$, where
$P_{\nu}(z)$ is the Legendre function of the first kind}
\author{Rados{\l}aw Szmytkowski \\*[3ex]
Atomic Physics Division,
Department of Atomic Physics and Luminescence, \\
Faculty of Applied Physics and Mathematics,
Gda{\'n}sk University of Technology, \\
Narutowicza 11/12, 80--233 Gda{\'n}sk, Poland \\
Email: radek@mif.pg.gda.pl}
\date{\today}
\maketitle
\begin{abstract}
We derive explicit expressions for the parameter derivatives 
$[\partial^{2}P_{\nu}(z)/\partial\nu^{2}]_{\nu=0}$ and
$[\partial^{3}P_{\nu}(z)/\partial\nu^{3}]_{\nu=0}$, where 
$P_{\nu}(z)$ is the Legendre function of the first kind. It is found 
that
\begin{displaymath}
\frac{\partial^{2}P_{\nu}(z)}{\partial\nu^{2}}\bigg|_{\nu=0}
=-2\Li_{2}\frac{1-z}{2},
\end{displaymath}
where $\Li_{2}z$ is the dilogarithm (this formula has been recently
arrived at by Schramkowski using \emph{Mathematica}), and that
\begin{displaymath}
\frac{\partial^{3}P_{\nu}(z)}{\partial\nu^{3}}\bigg|_{\nu=0}
=12\Li_{3}\frac{z+1}{2}-6\ln\frac{z+1}{2}\Li_{2}\frac{z+1}{2}
-\pi^{2}\ln\frac{z+1}{2}-12\zeta(3),
\end{displaymath}
where $\Li_{3}z$ is the polylogarithm of order 3 and $\zeta(s)$ is
the Riemann zeta function.
\end{abstract}
%
%
\section{Introduction}
\label{I}
\setcounter{equation}{0}
Some time ago, we carried our research on the first-order derivative
of the Legendre function of the first kind, $P_{\nu}(z)$, with
respect to its degree $\nu$ \cite{Szmy06,Szmy07}. Recently,
Schramkowski \cite{Schr12} has drawn our attention to the fact that
in the course of modeling of tidal hydrodynamics one has to
approximate the function $P_{\nu}(z)$ by a few first terms of its
Maclaurin series in $\nu$:
\begin{equation}
P_{\nu}(z)=\sum_{k=0}^{\infty}\frac{\nu^{k}}{k!}
\frac{\partial^{k}P_{\nu}(z)}{\partial\nu^{k}}\bigg|_{\nu=0},
\label{1.1}
\end{equation}
and that at least the terms with $k\leqslant3$ have to be retained if
results of the modeling are to be physically correct. Using
\emph{Mathematica\/}, Schramkowski has found that
\begin{equation}
\frac{\partial^{2}P_{\nu}(z)}{\partial\nu^{2}}\bigg|_{\nu=0}
=-2\Li_{2}\frac{1-z}{2},
\label{1.2}
\end{equation}
where $\Li_{2}z$ is the dilogarithm function \cite{Apos10}. In this
note, we provide a proof of the formula in Eq.\ (\ref{1.2}). In
addition, we derive a closed-form expression for the third-order
derivative $[\partial^{3}P_{\nu}(z)/\partial\nu^{3}]_{\nu=0}$.
\section{The derivative
$[\partial^{2}P_{\nu}(z)/\partial\nu^{2}]_{\nu=0}$}
\label{II}
\setcounter{equation}{0}
The Legendre function of the first kind, $P_{\nu}(z)$, obeys the
differential identity
\begin{equation}
\left[\frac{\mathrm{d}}{\mathrm{d}z}(1-z^{2})
\frac{\mathrm{d}}{\mathrm{d}z}+\nu(\nu+1)\right]P_{\nu}(z)=0.
\label{2.1}
\end{equation}
Differentiating Eq.\ (\ref{2.1}) twice with respect to $\nu$ and
rearranging terms gives
\begin{equation}
\left[\frac{\mathrm{d}}{\mathrm{d}z}(1-z^{2})
\frac{\mathrm{d}}{\mathrm{d}z}+\nu(\nu+1)\right]
\frac{\partial^{2}P_{\nu}(z)}{\partial\nu^{2}}
=-2P_{\nu}(z)-2(2\nu+1)
\frac{\partial P_{\nu}(z)}{\partial\nu}.
\label{2.2}
\end{equation}
One may consider Eq.\ (\ref{2.2}) as an inhomogeneous second-order
linear differential equation, a particular solution to which is
$\partial^{2}P_{\nu}(z)/\partial\nu^{2}$.

In the particular case of $\nu=0$, Eq.\ (\ref{2.2}) simplifies to
\begin{equation}
\frac{\mathrm{d}}{\mathrm{d}z}(1-z^{2})
\frac{\mathrm{d}}{\mathrm{d}z}
\frac{\partial^{2}P_{\nu}(z)}{\partial\nu^{2}}\bigg|_{\nu=0}
=-2P_{0}(z)-2\frac{\partial P_{\nu}(z)}{\partial\nu}\bigg|_{\nu=0}.
\label{2.3}
\end{equation}
Since it is known that
\begin{equation}
P_{0}(z)=1
\label{2.4}
\end{equation}
and \cite{Szmy06}
\begin{equation}
\frac{\partial P_{\nu}(z)}{\partial\nu}\bigg|_{\nu=0}
=\ln\frac{z+1}{2},
\label{2.5}
\end{equation}
Eq.\ (\ref{2.3}) may be rewritten as
\begin{equation}
\frac{\mathrm{d}}{\mathrm{d}z}(1-z^{2})
\frac{\mathrm{d}}{\mathrm{d}z}
\frac{\partial^{2}P_{\nu}(z)}{\partial\nu^{2}}\bigg|_{\nu=0}
=-2-2\ln\frac{z+1}{2}.
\label{2.6}
\end{equation}
Integrating both sides of Eq.\ (\ref{2.6}) with respect to $z$ yields
\begin{equation}
(1-z^{2})\frac{\mathrm{d}}{\mathrm{d}z}
\frac{\partial^{2}P_{\nu}(z)}{\partial\nu^{2}}\bigg|_{\nu=0}
=-2z-2(z+1)\left(\ln\frac{z+1}{2}-1\right)+C_{1},
\label{2.7}
\end{equation}
where $C_{1}$ is an integration constant, or equivalently as
\begin{equation}
(1-z^{2})\frac{\mathrm{d}}{\mathrm{d}z}
\frac{\partial^{2}P_{\nu}(z)}{\partial\nu^{2}}\bigg|_{\nu=0}
=-2(z+1)\ln\frac{z+1}{2}+C_{1}^{\prime},
\label{2.8}
\end{equation}
where $C_{1}^{\prime}=C_{1}+2$. Dividing both sides of Eq.\
(\ref{2.8}) by $1-z^{2}$ and integrating the resulting equation with
respect to $z$ furnishes the following expression for
$[\partial^{2}P_{\nu}(z)/\partial\nu^{2}]_{\nu=0}$:
\begin{equation}
\frac{\partial^{2}P_{\nu}(z)}{\partial\nu^{2}}\bigg|_{\nu=0}
=-\int\mathrm{d}z\:\frac{\displaystyle\ln\frac{z+1}{2}}
{\displaystyle\frac{1-z}{2}}
+\frac{C_{1}^{\prime}}{2}\ln\frac{z+1}{z-1}.
\label{2.9}
\end{equation}
However, it holds that
\begin{equation}
\int\mathrm{d}z\:\frac{\displaystyle\ln\frac{z+1}{2}}
{\displaystyle\frac{1-z}{2}}
=-2\int\mathrm{d}\left(\frac{1-z}{2}\right)\:
\frac{\displaystyle\ln\left(1-\frac{1-z}{2}\right)}
{\displaystyle\frac{1-z}{2}}=2\Li_{2}\frac{1-z}{2}+C_{2},
\label{2.10}
\end{equation}
where
\begin{equation}
\Li_{2}z=-\int_{0}^{z}\mathrm{d}t\:\frac{\ln(1-t)}{t}
\label{2.11}
\end{equation}
is the dilogarithm function \cite{Apos10} and $C_{2}$ is an
integration constant. Hence, we have
\begin{equation}
\frac{\partial^{2}P_{\nu}(z)}{\partial\nu^{2}}\bigg|_{\nu=0}
=-2\Li_{2}\frac{1-z}{2}+\frac{C_{1}^{\prime}}{2}\ln\frac{z+1}{z-1}
-C_{2}.
\label{2.12}
\end{equation}

To fix values of the constants $C_{1}^{\prime}$ and $C_{2}$, we use
the identity
\begin{equation}
P_{\nu}(1)=1,
\label{2.13}
\end{equation}
from which it follows that
\begin{equation}
\frac{\partial^{2}P_{\nu}(1)}{\partial\nu^{2}}\bigg|_{\nu=0}=0.
\label{2.14}
\end{equation}
The logarithmic function multiplying $C_{1}^{\prime}$ on the
right-hand side of Eq.\ (\ref{2.12}) is singular at $z=1$ and
consequently the necessary condition for Eq.\ (\ref{2.14}) holds is
\begin{equation}
C_{1}^{\prime}=0.
\label{2.15}
\end{equation}
Furthermore, we see that in view of the property
\begin{equation}
\Li_{2}0=0,
\label{2.16}
\end{equation}
it must hold that
\begin{equation}
C_{2}=0.
\label{2.17}
\end{equation}
This leads us to the final result
\begin{equation}
\frac{\partial^{2}P_{\nu}(z)}{\partial\nu^{2}}\bigg|_{\nu=0}
=-2\Li_{2}\frac{1-z}{2}.
\label{2.18}
\end{equation}
In this way, we have confirmed the finding of Schramkowski.
\section{The derivative
$[\partial^{3}P_{\nu}(z)/\partial\nu^{3}]_{\nu=0}$}
\label{III}
\setcounter{equation}{0}
The reasoning analogous to that presented in Sec.\ \ref{II} may be
carried out to evaluate the third-order derivative 
$[\partial^{3}P_{\nu}(z)/\partial\nu^{3}]_{\nu=0}$. Differentiating
Eq.\ (\ref{2.1}) triply with respect to $\nu$ yields
\begin{equation}
\left[\frac{\mathrm{d}}{\mathrm{d}z}(1-z^{2})
\frac{\mathrm{d}}{\mathrm{d}z}+\nu(\nu+1)\right]
\frac{\partial^{3}P_{\nu}(z)}{\partial\nu^{3}}
=-6\frac{\partial P_{\nu}(z)}{\partial\nu}
-3(2\nu+1)\frac{\partial^{2}P_{\nu}(z)}{\partial\nu^{2}},
\label{3.1}
\end{equation}
which for $\nu=0$ simplifies to 
\begin{equation}
\frac{\mathrm{d}}{\mathrm{d}z}(1-z^{2})
\frac{\mathrm{d}}{\mathrm{d}z}
\frac{\partial^{3}P_{\nu}(z)}{\partial\nu^{3}}\bigg|_{\nu=0}
=-6\frac{\partial P_{\nu}(z)}{\partial\nu}\bigg|_{\nu=0}
-3\frac{\partial^{2}P_{\nu}(z)}{\partial\nu^{2}}\bigg|_{\nu=0}.
\label{3.2}
\end{equation}
Use of Eqs.\ (\ref{2.5}) and (\ref{2.18}) transforms Eq.\ (\ref{3.2})
into the more explicit form
\begin{equation}
\frac{\mathrm{d}}{\mathrm{d}z}(1-z^{2})
\frac{\mathrm{d}}{\mathrm{d}z}
\frac{\partial^{3}P_{\nu}(z)}{\partial\nu^{3}}\bigg|_{\nu=0}
=-6\ln\frac{z+1}{2}+6\Li_{2}\frac{1-z}{2}.
\label{3.3}
\end{equation}
Integrating both sides of Eq.\ (\ref{3.3}) with respect to $z$ and
using
\begin{equation}
\int\mathrm{d}z\:\Li_{2}z=z\Li_{2}z+(z-1)\ln(1-z)-z+C_{1}
\label{3.4}
\end{equation}
(this results may be easily arrived at through integration by parts)
gives
\begin{equation}
(1-z^{2})\frac{\mathrm{d}}{\mathrm{d}z}
\frac{\partial^{3}P_{\nu}(z)}{\partial\nu^{3}}\bigg|_{\nu=0}
=6(z-1)\Li_{2}\frac{1-z}{2}+C_{1}^{\prime},
\label{3.5}
\end{equation}
hence, it follows that
\begin{equation}
\frac{\partial^{3}P_{\nu}(z)}{\partial\nu^{3}}\bigg|_{\nu=0}
=-6\int\mathrm{d}z\:\frac{\displaystyle\Li_{2}\frac{1-z}{2}}{z+1}
+\frac{C_{1}^{\prime}}{2}\ln\frac{z+1}{z-1}.
\label{3.6}
\end{equation}
Now, it is known \cite[Eq.\ (1.2.3.1)]{Prud03} that
\begin{equation}
\int\mathrm{d}z\:\frac{\Li_{2}z}{1-z}
=2\Li_{3}(1-z)-\ln(1-z)\Li_{2}z
-2\ln(1-z)\Li_{2}(1-z)-\ln z\ln^{2}(1-z)+C_{2},
\label{3.7}
\end{equation}
where $\Li_{3}z$ is the polylogarithm \cite{Apos10} of order 3.
Employing the Euler's identity \cite[p.\ 652]{Prud03}
\begin{equation}
\Li_{2}z+\Li_{2}(1-z)=\frac{\pi^{2}}{6}-\ln z\ln(1-z),
\label{3.8}
\end{equation}
Eq.\ (\ref{3.7}) may be simplified to
\begin{equation}
\int\mathrm{d}z\:\frac{\Li_{2}z}{1-z}=2\Li_{3}(1-z)
-\ln(1-z)\Li_{2}(1-z)-\frac{\pi^{2}}{6}\ln(1-z)+C_{2}.
\label{3.9}
\end{equation}
Evaluating the integral on the right-hand side of Eq.\ (\ref{3.6})
with the aid of the formula in Eq.\ (\ref{3.9}) leads to
\begin{equation}
\frac{\partial^{3}P_{\nu}(z)}{\partial\nu^{3}}\bigg|_{\nu=0}
=12\Li_{3}\frac{z+1}{2}-6\ln\frac{z+1}{2}\Li_{2}\frac{z+1}{2}
-\pi^{2}\ln\frac{z+1}{2}+\frac{C_{1}^{\prime}}{2}\ln\frac{z+1}{z-1}
+C_{2}^{\prime}.
\label{3.10}
\end{equation}
Since it holds that
\begin{equation}
\Li_{3}1=\sum_{k=1}^{\infty}\frac{1}{k^{3}}=\zeta(3),
\label{3.11}
\end{equation}
where $\zeta(s)$ is the Riemann zeta function, to ensure that the
condition
\begin{equation}
\frac{\partial^{3}P_{\nu}(1)}{\partial\nu^{3}}\bigg|_{\nu=0}=0
\label{3.12}
\end{equation}
(resulting from Eq.\ (\ref{2.13})) is satisfied, we must put
\begin{equation}
C_{1}^{\prime}=0
\label{3.13}
\end{equation}
and
\begin{equation}
C_{2}^{\prime}=-12\zeta(3).
\label{3.14}
\end{equation}
Hence, the final expression for 
$[\partial^{3}P_{\nu}(z)/\partial\nu^{3}]_{\nu=0}$ is
\begin{equation}
\frac{\partial^{3}P_{\nu}(z)}{\partial\nu^{3}}\bigg|_{\nu=0}
=12\Li_{3}\frac{z+1}{2}-6\ln\frac{z+1}{2}\Li_{2}\frac{z+1}{2}
-\pi^{2}\ln\frac{z+1}{2}-12\zeta(3).
\label{3.15}
\end{equation}
\section*{Acknowledgments}
I thank Dr.\ George P.\ Schramkowski for drawing my attention to the
problem and for communicating his finding concerning the explicit
form of the derivative
$[\partial^{2}P_{\nu}(z)/\partial\nu^{2}]_{\nu=0}$.
\end{document}